\documentclass[final,leqno]{siamltex}
\usepackage{amsfonts,amsmath,amssymb,graphicx,latexsym,makeidx,url}
\allowdisplaybreaks
\newtheorem{remark}[theorem]{Remark}
\newtheorem{algorithm}[theorem]{Algorithm}
\newcommand{\R}{\in \mathbb{R}}

\newcommand{\Rpp}{\R_{++}}
\newcommand{\C}{\in \mathbb{C}}
\newcommand{\Cpp}{\mathbb{C}_{++}}
\newcommand{\Zp}{\in \mathbb{Z}_+}
\newcommand{\Zpp}{\in \mathbb{Z}_{++}}
\newcommand{\ap}{\tilde}
\newcommand{\lam}{\lambda}
\newcommand{\Lam}{\Lambda}
\newcommand{\itv}[2]{\langle#1,#2\rangle}
\newcommand{\Itv}[2]{\left\langle#1,#2\right\rangle}
\newcommand{\dsfrac}[2]{{\displaystyle \frac{#1}{#2}}}
\newcommand{\norm}[1]{\|#1\|}

\newcommand{\normi}[1]{\|#1\|_\infty}
\newcommand{\normp}[1]{\|#1\|_p}
\newcommand{\Normp}[1]{\left\|#1\right\|_p}
\newcommand{\inv}[1]{#1^{-1}}
\newcommand{\RM}[2]{\R^{#1 \times #2}}
\newcommand{\RMp}[2]{\R_+^{#1 \times #2}}
\newcommand{\RV}[1]{\R^{#1}}
\newcommand{\RVp}[1]{\R_+^{#1}}
\newcommand{\M}[2]{\C^{#1 \times #2}}
\newcommand{\V}[1]{\C^{#1}}
\newcommand{\abs}[1]{|#1|}
\newcommand{\Abs}[1]{\left|#1\right|}
\newcommand{\Zm}{\mathbb{Z}_-}
\newcommand{\Od}{\mathcal{O}}
\newcommand{\bm}[1]{\mbox{\boldmath $#1$}}
\newcommand{\om}[1]{{\rm 1}\hspace{-0.25em}{\rm l}_{#1}^{\rm M}}
\newcommand{\ov}[1]{{\rm 1}\hspace{-0.25em}{\rm l}_{#1}^{\rm v}}
\newcommand{\rr}{r_{{\rm r}}}
\newcommand{\rc}{r_{{\rm c}}}
\newcommand{\su}[1]{^{(#1)}}
\newcommand{\IM}[2]{\IC^{#1 \times #2}}
\newcommand{\IC}{\in \mathbb{IC}}
\title{Verified computation of matrix gamma function 
\thanks{This work was partially supported by JSPS KAKENHI Grant Number JP16K05270.}}
\author{Shinya Miyajima\thanks{Faculty of Science and Engineering, Iwate University,  020-8551, Japan ({\tt miyajima@iwate-u.ac.jp}).}}
\begin{document}
\maketitle
\begin{abstract}
Two numerical algorithms are proposed for computing an interval matrix containing the matrix gamma function. 
In 2014, the author presented algorithms for enclosing all the eigenvalues and basis of invariant subspaces of $A \M{n}{n}$. 
As byproducts of these algorithms, we can obtain interval matrices containing small matrices whose spectrums are included in that of $A$. 
In this paper, we interpret the interval matrices containing the basis and small matrices as a result of verified block diagonalization (VBD), 
and establish a new framework for enclosing matrix functions using the VBD. 
To achieve enclosure for the gamma function of the small matrices, we derive computable perturbation bounds.  
We can apply these bounds if input matrices satisfy conditions. 
We incorporate matrix argument reductions (ARs) to force the input matrices to satisfy the conditions, and develop theories for accelerating the ARs. 
The first algorithm uses the VBD  based on a numerical spectral decomposition, and involves only cubic complexity under an assumption. 
The second algorithm adopts the VBD based on a numerical Jordan decomposition, and is applicable even for defective matrices. 
Numerical results show efficiency and robustness of the algorithms. 
\end{abstract}
%
\pagestyle{myheadings}
\thispagestyle{plain}
\markboth{SHINYA MIYAJIMA}{VERIFIED COMPUTATION OF MATRIX GAMMA FUNCTION}
\section{Introduction}\label{sec:I}
For $z \C$ with positive real part, the gamma function is defined by 
$$
\Gamma(z) := \int_0^\infty e^{-t}t^{z-1}dt,
$$
and otherwise by analytic continuation. 
It is well known that $\Gamma(z)$ is analytic everywhere in $\mathbb{C}$, with the exception of non-positive integer numbers $\Zm$. 
Therefore, the general theory of primary matrix function \cite{H08} ensures that the matrix gamma function $\Gamma(A)$ is well defined for $A \M{n}{n}$ having no eigenvalues on $\Zm$. 
If all eigenvalues of $A$ have positive real parts, then we have the representation 
\begin{equation}\label{eq:Gamma}
\Gamma(A) = \int_0^\infty e^{-t}t^{A-I_n}dt, 
\end{equation}
where $t^{A-I_n} := e^{(A-I_n)\log(t)}$ and $I_n$ denotes the $n \times n$ identity matrix. 

The function $\Gamma(A)$ has connections with other special functions, which play an important role in solving certain matrix differential equations \cite{CS}. 
Two of these special functions are the matrix beta and Bessel functions. 
In \cite{CS}, mathematical properties of $\Gamma(A)$ are elegantly clarified, and fast and accurate algorithms for computing $\Gamma(A)$ are proposed. 

The work presented in this paper addresses the problem of verified computations for $\Gamma(A)$, 
specifically, numerically computing interval matrices which are guaranteed to contain $\Gamma(A)$. 
To the author's best knowledge, a verification algorithm designed specifically for $\Gamma(A)$ does not yet appear in the literature. 
A possible method is to use the VERSOFT \cite{R} routine {\tt vermatfun}. 
This routine is applicable not only to the matrix gamma function but also to other matrix functions, 
and computes the interval matrices by enclosing all the eigenvalues and eigenvectors of $A$ via the INTLAB  \cite{R99} routine {\tt verifyeig}. 
This routine fails when $A$ is defective or close to defective, and requires $\Od(n^4)$ operations. 

The purpose of this paper is to propose two verification algorithms for $\Gamma(A)$. 
In \cite{M14}, algorithms for enclosing all the eigenvalues and basis of invariant subspaces of $A$ are presented. 
As byproducts of these algorithms, we can obtain interval matrices containing small matrices whose spectrums are included in that of $A$. 
In this paper, we interpret the interval matrices containing the basis and small matrices as a result of verified block diagonalization (VBD), 
and establish a new framework for enclosing matrix functions using the VBD. 
To achieve enclosure for the gamma function of the small matrices, we derive computable perturbation bounds. 
Here, the word ``computable'' means that we can numerically obtain a rigorous upper bound which takes rounding and truncation errors into account. 
We can find a perturbation bound for $\Gamma(A)$ also in \cite{CS}. 
On the other hand, the bound in \cite{CS} is not a computable one. 
We can apply the derived perturbation bounds if input matrices satisfy conditions. 
We incorporate matrix argument reductions (ARs) to force the input matrices to satisfy the conditions, and develop theories for accelerating the ARs. 
The first algorithm uses the VBD based on a numerical spectral decomposition (NSD), and involves only $\Od(n^3)$ operations under an assumption. 
The second algorithm adopts the VBD based on a numerical Jordan decomposition (NJD), and is applicable even when $A$ is defective. 
We present a theory for verifying that $A$ has no eigenvalues on $\Zm$. 
By the aid of this theory, these algorithms do not assume but prove that $A$ has no eigenvalues on $\Zm$. 
The first and second algorithms require intervals containing $\Gamma\su{0}(z)/0!,\dots,\Gamma\su{\ell}(z)/\ell!$, where $\ell$ is a non-negative integer and $z \C$. 
To the author's best knowledge, an algorithm for computing such intervals is not available in literature, 
whereas there are well-established algorithms \cite{K,R14,YOO} for computing intervals containing {\em real} scalar gamma functions. 
We thus present a way for computing such intervals, which is based on the Spouge approximation \cite{S}. 
Although this way may be a slight modification of the Spouge method, the proposed algorithms are the first ones which apply the VBD to computation of an interval containing a matrix function. 
One may consider that the VBD is a direct application of the algorithms in \cite{M14}. 
However, the established framework enables us to enclose not only $\Gamma(A)$ but also other matrix functions (see Section~\ref{sec:C}). 
Moreover, utilizing the VBD as a means to enclose a matrix function, verifying that $A$ has no eigenvalues on $\Zm$, deriving the computable perturbation bounds, 
and accelerating of the ARs are the first attempts and not obvious. 

The author has been proposed many verification algorithms for matrix functions (e.g., \cite{M18rt,M18ex,M19ln,M19w}). 
However, the idea in this paper does not overlap with those in the previous papers. 
This is because most of the previous algorithms are based on matrix equations, whereas the algorithms in this paper are not. 
Although the algorithms in \cite{M18ex} are not based on matrix equations and also utilize the NSD or NJD, these algorithms do not use the VBD, which is the key idea in this paper.  

This paper is organized as follows: 
Section~\ref{sec:P} introduces notation and theories used in this paper. 
Section~\ref{sec:GD} presents a way for computing the intervals containing $\Gamma\su{0}(z)/0!,\dots,\Gamma\su{\ell}(z)/\ell!$. 
Sections~\ref{sec:S} and \ref{sec:J} propose the first and second algorithms, respectively. 
Section~\ref{sec:N} reports numerical results. 
Section~\ref{sec:C} finally summarizes the results in this paper and highlights possible extension and future work. 
\section{Preliminaries}\label{sec:P}
For $M \M{n}{n}$, let $M_{ij}$, $M_{:j}$, $\rho(M)$ and $\mu(M)$ be the $(i,j)$ element, $j$-th column, spectral radius and spectrum of $M$, respectively, and $\abs{M} := (\abs{M_{ij}})$. 
For $v \V{n}$, denote the $i$-th element of $v$ by $v_i$. 
For $M,N \RM{m}{n}$, the inequality $M \le N$ means $M_{ij} \le N_{ij}$, $\forall i,j$. 
Let ${\sf i} := \sqrt{-1}$, $\mathbb{Z}_+ := \{z \in \mathbb{Z} :  z \ge 0\}$, $\mathbb{Z}_{++} := \{z \in \mathbb{Z} :  z > 0\}$, $\Zm := \{z \in \mathbb{Z} :  z \le 0\}$, 
$\mathbb{R}_{+} := [0,\infty)$, $\mathbb{R}_{++} := (0,\infty)$, $\Cpp := \{z \C : {\rm Re}(z) > 0\}$, $\mathbb{R}_+^n := \{v \RV{n} : v \ge 0\}$, 
and $\mathbb{R}_+^{m \times n} := \{M \RM{m}{n} : M \ge 0\}$. 
Let also $\mathbb{IC}$ and $\mathbb{IC}^{m \times n}$ be the sets of all complex interval scalars and $m \times n$ matrices, respectively. 
For $C \M{m}{n}$ and $R \RMp{m}{n}$, denote the interval matrix whose midpoint and radius are $C$ and $R$, respectively, by $\itv{C}{R}$. 
Suppose any matrices contained in $\bm{M} \IM{n}{n}$ is nonsingular. 
Then, $\inv{\bm{M}}$ denotes an interval matrix including $\{\inv{M} : M \in \bm{M}\}$. 
Expressions containing intervals mean results of interval arithmetic. 
Let $A,B \M{n}{n}$ and $R,S \RMp{n}{n}$. 
In Sections~\ref{sec:S} and \ref{sec:J}, we will use the following property of interval arithmetic (see \cite{A}, e.g.): 
\begin{equation}\label{eq:A}
\itv{A}{R}\itv{B}{S} \subseteq \itv{AB}{\abs{A}S + R\abs{B} + RS}. 
\end{equation}
For $\alpha \R$, let $\lceil\alpha\rceil$ and $\lfloor\alpha\rfloor$ denote the ceiling and floor functions, respectively. 
In Sections~\ref{sec:S} and \ref{sec:J}, we will use the {\em incomplete} gamma function 
$$
\gamma(\alpha) := \int_0^1e^{-t}t^{\alpha-1}dt, \quad \mbox{where} \quad \alpha \Rpp.  
$$
For $z \C$, let $\log(z)$ be the principal branch of the logarithm. 
Define
$$
\ov{n} := \left[\begin{array}{c}
1 \\
\vdots \\
1\\
\end{array}\right] \RV{n}, \  
\om{n} := \left[\begin{array}{ccc}
1 & \cdots & 1 \\
\vdots &  & \vdots \\
1 & \cdots & 1 \\
\end{array}\right] \RM{n}{n} \  \mbox{and} \  N_j := \left[\begin{array}{cccc}
0 & 1 &  & \\
  & 0 & \ddots & \\
    & & \ddots & 1 \\
  & & & 0
\end{array}\right] \RM{n_j}{n_j}. 
$$

In Section~\ref{sec:GD}, we will use the Spouge approximation \cite{S} and its error bound, which are summarized in Lemma~\ref{lm:S}. 
\begin{lemma}[Spouge~\cite{S}]\label{lm:S}
Let $a \Rpp$ and $z,w \C$. 
Define $c_0 := 1$, 
\begin{eqnarray*}
c_k &:=& \frac{1}{\sqrt{2\pi}}\frac{(-1)^{k-1}}{(k-1)!}(-k + a)^{k - 1/2}e^{-k+a}, \quad  k = 1,2,\dots, \lceil a \rceil - 1,  \\
H(z) &:=& c_0 + \sum_{k=1}^{\lceil a \rceil - 1}\frac{c_k}{z-1+k}, \quad K(z) := \sqrt{2\pi}(z - 1 + a)^{z-1/2}e^{-(z-1+a)}, \\
G(w) &:=& \frac{\sqrt{2\pi}(-w-a)^{-w-1/2}e^{w+a}}{(-w-1)!}, \\
\epsilon(z) &:=& \frac{{\sf i}}{2\pi}\int_0^\infty\left( \frac{G(-a - {\sf i}v)}{e^{2\pi(v-{\sf i}a)}-1} + \frac{G(-a + {\sf i}v)}{e^{2\pi(v+{\sf i}a)}-1}\right)\frac{dv}{-a + {\sf i}v - z + 1}.
\end{eqnarray*}
Assume $a \ge 3$ and ${\rm Re}(z - 1 + a) > 0$. 
Then, 
\begin{description}
\item[{\rm (a)}] $\Gamma(z) = K(z)(H(z) + \epsilon(z))$;
\item[{\rm (b)}] for $m \Zp$, the $m$-th derivative of the error term $\epsilon(z)$ is bounded by 
$$
\abs{\epsilon\su{m}(z)} \le \frac{m!C_a}{({\rm Re}(z-1+a))^{m+1}}, \  \mbox{where} \ C_a := \frac{\sqrt{2/\pi}}{(a-1)!}\int_0^\infty \frac{v^{a-1/2}}{\abs{e^{2\pi v} - e^{2\pi {\sf i}a}}}dv;
$$ 
\item[{\rm (c)}] $C_a < \sqrt{ae/\pi}(2\pi)^{-(a+1/2)}$. 
\end{description}
\end{lemma}
\begin{remark}
We can obtain Lemma~\ref{lm:S} (c) from \cite[Proof of Theorem~1.3.1]{S}. 
\end{remark}
From Lemma~\ref{lm:S} (b) and (c), we immediately obtain Corollary~\ref{cr:S}. 
\begin{corollary}\label{cr:S} 
Let $m$, $a$, $z$ and $\epsilon(z)$ be as in Lemma~\ref{lm:S}. 
If $a \ge 3$ and ${\rm Re}(z - 1 + a) > 0$, then $\abs{\epsilon\su{m}(z)} < \xi_m(z)$, where 
$$
\xi_m(z) := \frac{m!\sqrt{ae/\pi}}{({\rm Re}(z-1+a))^{m+1}(2\pi)^{a+1/2}}.
$$
\end{corollary}

In Sections~\ref{sec:S} and \ref{sec:J}, we will use the following properties of matrix functions: 
\begin{lemma}[e.g., Higham \cite{H08}]\label{lm:H} 
Let $A,X,Y \M{n}{n}$ and $\varphi$ be defined on the spectrum of $A$. 
Then, 
\begin{description}
\item[{\rm (a)}] if $X$ is nonsingular, then $\varphi(XA\inv{X}) = X\varphi(A)\inv{X}$;
\item[{\rm (b)}] if $A = {\rm diag}(A_1,\dots,A_p)$ is block diagonal, then $\varphi(A) = {\rm diag}(\varphi(A_1),\dots,\varphi(A_p))$; 
\item[{\rm (c)}] if $XY = YX$, then $e^{X+Y} = e^Xe^Y = e^Ye^X$; 
\item[{\rm (d)}] $\norm{e^X - e^Y} \le \norm{X - Y}e^{\max(\norm{X},\norm{Y})}$ for any consistent norm. 
\end{description}
\end{lemma}

We cite Lemma~\ref{lm:CS} as a theoretical basis for the ARs in Sections~\ref{sec:ARS} and \ref{sec:ARJ}. 
\begin{lemma}[Cardoso and Sadeghi \cite{CS}]\label{lm:CS} 
Let $A \M{n}{n}$ have no eigenvalues on $\Zm$. 
Then, $\Gamma(A + I_n) = A\Gamma(A)$. 
\end{lemma}

Let $\alpha \Rpp$. 
In Sections~\ref{sec:EBS} and \ref{sec:EBJ}, we will estimate an upper bound for $-\gamma\su{1}(\alpha)$. 
To this end, we present Lemma~\ref{lm:IG}. 
\begin{lemma}\label{lm:IG}  
Let $\alpha \Rpp$ and $\omega(\alpha) := \dsfrac{2\alpha+1}{\alpha^2(\alpha+1)^2} + \dsfrac{\cosh(1)-1}{(\alpha + 2)^2}$. 
Then, $-\gamma\su{1}(\alpha) < \omega(\alpha)$. 
\end{lemma}
\proof From $\gamma\su{1}(\alpha) = \int_0^1e^{-t}t^{\alpha-1}\log(t)dt$ and integration by parts, we obtain
\begin{eqnarray}
-\gamma\su{1}(\alpha) 
&=& -\left[\log(t)\sum_{i=0}^\infty\frac{(-1)^it^{\alpha+i}}{i!(\alpha+i)}\right]_{t=0}^{t=1} + \int_0^1\left(\sum_{i=0}^\infty\frac{(-1)^it^{\alpha+i-1}}{i!(\alpha+i)}\right)dt \nonumber \\
&=& \sum_{i=0}^\infty\frac{(-1)^i}{i!(\alpha+i)^2} \nonumber \\
&=& \frac{1}{\alpha^2} - \frac{1}{(\alpha+1)^2} + \sum_{i=1}^\infty\left(\frac{1}{(2i)!(\alpha+2i)^2} - \frac{1}{(2i+1)!(\alpha+2i+1)^2}\right). \label{eq:IG}
\end{eqnarray}
For $i = 1,2,\dots$, it follows that 
\begin{eqnarray*}
\frac{1}{(2i)!(\alpha+2i)^2} - \frac{1}{(2i+1)!(\alpha+2i+1)^2} &=& \frac{2i}{(2i+1)!(\alpha+2i+1)^2} + \frac{2(\alpha+2i)+1}{(2i)!(\alpha+2i)^2(\alpha+2i+1)^2} \\
&<& \frac{1}{(2i)!(\alpha+3)^2} + \frac{2(\alpha+2)+1}{(2i)!(\alpha+2)^2(\alpha+3)^2} \\
&=& \frac{1}{(2i)!(\alpha+2)^2}. 
\end{eqnarray*}
This and (\ref{eq:IG}) prove the inequality. \quad \endproof
\section{Enclosing $\bm{\Gamma\su{0}(z)/0!,\dots,\Gamma\su{\ell}(z)/\ell!}$}\label{sec:GD}
As mentioned in Section~\ref{sec:I}, we need to compute intervals containing $\Gamma\su{0}(z)/0!,\dots,\Gamma\su{\ell}(z)/\ell!$ for $z \C$ and $\ell \Zp$. 
To this end, we use Lemma~\ref{lm:S} and Corollary~\ref{cr:S}. 
\begin{remark}
There are many other methods for computing an approximation of $\Gamma(z)$ (see \cite{CS}, e.g.). 
By exploiting these methods, computing an interval containing $\Gamma(z)$ seems to be possible. 
On the other hand, error bounds regarding to the {\em derivatives} of $\Gamma(z)$ are explicitly written in \cite{S}. 
Therefore, the Spouge method is useful for our purpose. 
\end{remark}

Let $a$, $c_k$, $H(z)$, $K(z)$ and $\epsilon(z)$ be as in Lemma~\ref{lm:S}, and $\xi_m(z)$ be as in Corollary~\ref{cr:S}. 
Suppose $a \ge 3$ and ${\rm Re}(z - 1 + a) > 0$. 
From Lemma~\ref{lm:S}, Corollary~\ref{cr:S}, and the Leibniz rule, for $m = 0,\dots,\ell$, we have
\begin{eqnarray}
\frac{\Gamma\su{m}(z)}{m!} &=& \sum_{k=0}^m\frac{K\su{k}(z)}{k!}\left(\frac{H^{(m-k)}(z)}{(m-k)!} + \frac{\epsilon^{(m-k)}(z)}{(m-k)!}\right) \nonumber \\
&\in& \sum_{k=0}^m\frac{K\su{k}(z)}{k!}\Itv{\frac{H^{(m-k)}(z)}{(m-k)!}}{\frac{\abs{\epsilon^{(m-k)}(z)}}{(m-k)!}}  \nonumber \\
&\subseteq& \sum_{k=0}^m\frac{K\su{k}(z)}{k!}\Itv{\frac{H^{(m-k)}(z)}{(m-k)!}}{\frac{\xi_{m-k}(z)}{(m-k)!}}. \label{eq:EncGam} 
\end{eqnarray}
We thus enclose $H\su{k}(z)/k!$ and $K\su{k}(z)/k!$ for $k = 0,\dots,\ell$. 
For large $k$, on the other hand, explicit representations for $K\su{k}(z)/k!$ seems to be complicated. 
For enclosing $K\su{k}(z)/k!$ without using the explicit representations, we propose the following way: 
Let $P(z) := \log(z - 1 + a) - (a-1/2)/(z - 1 +a)$. 
Then, $K\su{1}(z) = K(z)P(z)$, so that 
\begin{equation}\label{eq:K}
\frac{K\su{k+1}(z)}{(k+1)!} = \frac{(K(z)P(z))\su{k}}{(k+1)!} = \frac{1}{k+1}\sum_{j=0}^k\frac{K\su{j}(z)}{j!}\frac{P\su{k-j}(z)}{(k-j)!}, \quad k = 0,\dots,\ell-1. 
\end{equation}
Hence, we can enclose $K^{(k+1)}(z)/(k+1)!$ if enclosures for $K\su{0}(z)/0!,\dots,K\su{k}(z)/k!$ have already been obtained. 
Observe that we can easily write down $H\su{j}(z)/j!$ and $P\su{j}(z)/j!$ explicitly. 
For $j = 1,\dots,\ell$, in fact, 
\begin{eqnarray}
\frac{H\su{j}(z)}{j!} &=& \sum_{k=1}^{\lceil a \rceil - 1}\frac{(-1)^jc_k}{(z - 1 + k)^{j+1}},  \label{eq:H} \\
\frac{P\su{j}(z)}{j!} &=& \frac{(-1)^{j-1}}{j(z - 1 + a)^j} + \frac{(-1)^{j-1}(a-1/2)}{(z - 1 + a)^{j+1}}.  \label{eq:P}
\end{eqnarray}
We summarize our approach in Algorithm~\ref{alg:G}. 
\begin{algorithm}\label{alg:G}
Let $a \ge 3$ be given and $\ell \Zp$. 
Assume ${\rm Re}(z - 1 + a) > 0$. 
This algorithm computes intervals containing $\Gamma\su{0}(z)/0!,\dots,\Gamma\su{\ell}(z)/\ell!$. 
\begin{description}
\item[Step 1.] Enclose $H\su{j}(z)/j!$ and $P\su{j}(z)/j!$ for $j = 0,\dots,\ell$ based on (\ref{eq:H}) and (\ref{eq:P}), respectively. 
\item[Step 2.] Compute intervals including $K\su{j}(z)/j!$ for $j = 0,\dots,\ell$ based on (\ref{eq:K}). 
\item[Step 3.] Enclose $\Gamma\su{j}(z)/j!$ for $j = 0,\dots,\ell$ based on (\ref{eq:EncGam}). 
\end{description}
\end{algorithm}
Step~1 involves $\Od(\lceil a \rceil\ell)$ operations. 
Steps~2 and 3 require $\Od(\ell^2)$ operations. 
Therefore, Algorithm~\ref{alg:G} involves $\Od(\lceil a \rceil\ell + \ell^2)$ operations. 

For executing Algorithm~\ref{alg:G}, we need to determine $a$. 
From the assumption in Corollary~\ref{cr:S}, we focus on the case $a \ge 3$. 
If we take $a$ too small, then $\xi_0(z)$ does not become small. 
If we take $a$ too large, on the other hand, many interval arithmetics are required for computing an interval containing $H(z)$, 
which causes enlargement of the radius of the interval. 
If we take $a$ in the form of $a = b + 1/2$, where $b \Zp$ is not too large, then rounding errors do not occur in the floating point computations of $a + 1/2$ and $a - 1/2$. 
Based on these observations, we propose incrementing $a$ by one from 7/2, and terminating the increment when the radius exceeds $\xi_0(z)$. 
We summarize this strategy in Algorithm~\ref{alg:a}. 
\begin{algorithm}\label{alg:a}
Assume ${\rm Re}(z) > -5/2$. 
This algorithm determines $a$ in Algorithm~\ref{alg:G}. 
\begin{description}
\item[Step 1.] Initialize $a$ as $a = 7/2$. 
\item[Step 2.] Compute intervals containing $c_1,\dots,c_{\lceil a\rceil-1}$ and $H(z)$. 
\item[Step 3.] If the radius of the interval containing $H(z)$ exceeds $\xi_0(z)$, then output the current $a$ and terminate. 
Otherwise, go to Step 4. 
\item[Step 4.] Update $a$ such that $a = a + 1$ and go back to Step~2. 
\end{description}
\end{algorithm}
Note that $c_1,\dots,c_{\lceil a\rceil-1}$ are computed whenever $a$ is incremented. 
Algorithm~\ref{alg:a} thus requires $\Od(\lceil a\rceil)$ operations per iteration. 

By slightly modifying Algorithms~\ref{alg:G} and \ref{alg:a}, we can compute intervals containing $\{\Gamma\su{j}(z)/j! : z \in \bm{z}\}$ for $j = 0,\dots,\ell$, 
where $\bm{z} =: \itv{c}{r} \IC$ satisfies ${\rm Re}(c) - r > -5/2$. 
To be specific, by replacing $z$ and $\xi_m(z)$ in these algorithms by $\bm{z}$ and  
$$
\frac{m!\sqrt{ae/\pi}}{({\rm Re}(c)-r-1+a)^{m+1}(2\pi)^{a+1/2}},
$$ 
respectively, we can obtain such intervals. 
\section{Algorithm based on the NSD}\label{sec:S}
We develop our algorithm in some steps. 
Section~\ref{sec:BDS} introduces the VBD based on the NSD, and framework using the VBD. 
Section~\ref{sec:VS} develops the theory for verifying $\mu(A) \cap \mathbb{Z}_- = \emptyset$. 
Section~\ref{sec:EBS} establishes the computable perturbation bound for enclosing the gamma function of a diagonal block. 
Section~\ref{sec:ARS} explains the ARs, and presents the theory for its acceleration. 
Section~\ref{sec:OAS} proposes the overall algorithm. 
\subsection{The VBD based on the NSD}\label{sec:BDS}  
Assume as a result of the NSD of $A$, 
we have $\Lam,X \M{n}{n}$ with $\Lam = {\rm diag}(\lam_1,\dots,\lam_n)$ such that $AX \approx X\Lam$. 
By executing column permutation if necessary, let $\{\lam_{i_1\su{j}},\dots,\lam_{i_{p_j}\su{j}}\}$, $j = 1,\dots,q$ be sets of clusters, where  
$i_1\su{1},\dots,i_{p_1}\su{1},\dots,i_1\su{q},\dots,i_{p_q}\su{q} \Zpp$ satisfy 
$1 = i_1\su{1} < \cdots < i_{p_1}\su{1} < \cdots < i_1\su{q} < \cdots < i_{p_q}\su{q} = n$ and $p_1 + \cdots + p_q = n$. 
Note that the case where $\lam_j$ is isolated from the others is included in the case $p_j = 1$. 
Let also $X_j := [X_{:i_1\su{j}},\dots,X_{:i_{p_j}\su{j}}]$ for $j = 1,\dots,q$, and $W_j \M{n}{p_j}$ and $P_j \M{p_j}{p_j}$ satisfy $AW_j = W_jP_j$. 
Observe $\mu(P_j) \subseteq \mu(A)$. 
Then, \cite[Algorithm~1]{M14} gives $\Delta_j \RMp{n}{p_j}$ and $\itv{\widehat{\lam}_j}{\varrho_j} \IC$ such that 
$\itv{X_j}{\Delta_j} \ni W_j$ and $\itv{\widehat{\lam}_j}{\varrho_j} \supseteq \mu(P_j)$ with $X_j$ and $\lam_{i_1\su{j}},\dots,\lam_{i_{p_j}\su{j}}$ being inputs, for $j = 1,\dots,q$. 
As byproducts of this algorithm, actually, we can obtain $R_j \RMp{p_j}{p_j}$ such that $\itv{\widehat{\lam}_jI_{p_j}}{R_j} \ni P_j$. 
Let $W := [W_1,\dots,W_q] \M{n}{n}$, $\Delta := [\Delta_1,\dots,\Delta_q] \RMp{n}{n}$ and $\bm{W} := \itv{X}{\Delta} \IM{n}{n}$. 
Then, $W \in \bm{W}$ and 
$$
AW = [AW_1,\dots,AW_q] = [W_1P_1,\dots,W_qP_q] = W{\rm diag}(P_1,\dots,P_q).
$$
We can verify nonsingularity of any matrix contained in $\bm{W}$ by executing a known algorithm (e.g., the INTLAB routine {\tt verifylss}). 
If the verification is succeeded, then $W$ is also nonsingular, so that $A = W{\rm diag}(P_1,\dots,P_q)\inv{W}$. 
Thus, $\bm{W}$ and ${\rm diag}(\itv{\widehat{\lam}_1I_{p_1}}{R_1},\dots,\itv{\widehat{\lam}_qI_{p_q}}{R_q})$ can be regarded as the result of the VBD. 
We establish the new framework for enclosing matrix functions based on the VBD. 
Although this paper treats $\Gamma(A)$ only, this framework enables us to enclose other matrix functions (see Section~\ref{sec:C}).  
From Lemma~\ref{lm:H} (a) and (b), we have $\Gamma(A) = W{\rm diag}(\Gamma(P_1),\dots,\Gamma(P_q))\inv{W}$, 
so that the problem of enclosing $\Gamma(A)$ can be reduced to that of enclosing $\Gamma(P_1),\dots,\Gamma(P_q)$. 
\subsection{Verification of $\bm{\mu(A) \cap \mathbb{Z}_- = \emptyset}$}\label{sec:VS}  
As another result of \cite[Algorithm~1]{M14}, we can obtain $r \RVp{n}$ such that $\mu(A) \subseteq \bigcup_{i=1}^n\itv{\lam_i}{r_i}$. 
We formulate and prove Theorem~\ref{th:VS} for verifying $\mu(A) \cap \mathbb{Z}_- = \emptyset$ using $\lam_i$ and $r$. 
\begin{theorem}\label{th:VS}  
Let $\lam \V{n}$ and $r \RVp{n}$ satisfy $\mu(A) \subseteq \bigcup_{i=1}^n\itv{\lam_i}{r_i}$. 
Define $f \V{n}$ by 
$$
f_i := \max({\rm Re}(\lam_i), \lfloor{\rm Re}(\lam_i)\rfloor - {\rm Re}(\lam_i), {\rm Re}(\lam_i) - \lceil{\rm Re}(\lam_i)\rceil) + {\rm Im}(\lam_i){\sf i}, \quad i = 1,\dots,n. 
$$ 
If $\min_i(\abs{f_i} - r_i) > 0$, then $\mu(A) \cap \mathbb{Z}_- = \emptyset$. 
\end{theorem}
\proof If $\bigcup_{i=1}^n\itv{\lam_i}{r_i} \cap \mathbb{Z}_- = \emptyset$, then $\mu(A) \cap \mathbb{Z}_- = \emptyset$. 
We thus prove $\itv{\lam_i}{r_i} \cap \mathbb{Z}_- = \emptyset$ for each $i$ by considering the cases of ${\rm Re}(\lam_i) \ge 0$ and ${\rm Re}(\lam_i) < 0$ separately. 

Consider first the case where ${\rm Re}(\lam_i) \ge 0$. 
Then, $\min_{c \in \Zm}\abs{\lam_i - c} = \abs{\lam_i}$, so that $\itv{\lam_i}{r_i} \cap \mathbb{Z}_- = \emptyset$ follows if $\abs{\lam_i} - r_i > 0$. 
Since $\lfloor{\rm Re}(\lam_i)\rfloor - {\rm Re}(\lam_i) \le 0$ and ${\rm Re}(\lam_i) - \lceil{\rm Re}(\lam_i)\rceil \le 0$, we have $f_i = \lam_i$. 
Therefore, $\abs{f_i} - r_i > 0$ is equivalent to $\abs{\lam_i} - r_i > 0$. 
Hence, $\itv{\lam_i}{r_i} \cap \mathbb{Z}_- = \emptyset$ if $\abs{f_i} - r_i > 0$. 

Consider next the case where ${\rm Re}(\lam_i) < 0$. 
Then, 
$$
\min_{c \in \Zm}\abs{\lam_i - c} = \min(\abs{\lfloor{\rm Re}(\lam_i)\rfloor - \lam_i},\abs{\lam_i - \lceil{\rm Re}(\lam_i)\rceil}) = : g_i, 
$$ 
so that  $\itv{\lam_i}{r_i} \cap \mathbb{Z}_- = \emptyset$ follows if $g_i - r_i > 0$. 
If ${\rm Re}(\lam_i) \in (-1,0)$, then $\lceil{\rm Re}(\lam_i)\rceil = 0$ gives ${\rm Re}(\lam_i) = {\rm Re}(\lam_i) - \lceil{\rm Re}(\lam_i)\rceil$. 
If ${\rm Re}(\lam_i) \le -1$, on the other hand, then $\lfloor{\rm Re}(\lam_i)\rfloor - {\rm Re}(\lam_i) \in (-1,0]$ and ${\rm Re}(\lam_i) - \lceil{\rm Re}(\lam_i)\rceil \in (-1,0]$ yield 
${\rm Re}(\lam_i) < \lfloor{\rm Re}(\lam_i)\rfloor - {\rm Re}(\lam_i)$ and ${\rm Re}(\lam_i) < {\rm Re}(\lam_i) - \lceil{\rm Re}(\lam_i)\rceil$. 
Therefore, $f_i$ can be written as 
$$
f_i = \max(\lfloor{\rm Re}(\lam_i)\rfloor - {\rm Re}(\lam_i), {\rm Re}(\lam_i) - \lceil{\rm Re}(\lam_i)\rceil) + {\rm Im}(\lam_i){\sf i}. 
$$ 
If $(\lfloor{\rm Re}(\lam_i)\rfloor + \lceil{\rm Re}(\lam_i)\rceil)/2 \le {\rm Re}(\lam_i)$, 
then $g_i = \abs{\lam_i - \lceil{\rm Re}(\lam_i)\rceil}$ and $\lfloor{\rm Re}(\lam_i)\rfloor - {\rm Re}(\lam_i) \le {\rm Re}(\lam_i) - \lceil{\rm Re}(\lam_i)\rceil$. 
Hence, $\abs{f_i} = \abs{{\rm Re}(\lam_i) - \lceil{\rm Re}(\lam_i)\rceil + {\rm Im}(\lam_i){\sf i}} = g_i$. 
If $(\lfloor{\rm Re}(\lam_i)\rfloor + \lceil{\rm Re}(\lam_i)\rceil)/2 > {\rm Re}(\lam_i)$, on the other hand, 
then $g_i = \abs{\lfloor{\rm Re}(\lam_i)\rfloor - \lam_i}$ and $\lfloor{\rm Re}(\lam_i)\rfloor - {\rm Re}(\lam_i) > {\rm Re}(\lam_i) - \lceil{\rm Re}(\lam_i)\rceil$. 
Thus, $\abs{f_i} = \abs{\lfloor{\rm Re}(\lam_i)\rfloor - {\rm Re}(\lam_i) + {\rm Im}(\lam_i){\sf i}} = g_i$. 
Therefore, if $\abs{f_i} - r_i > 0$, then $g_i - r_i > 0$, so that $\itv{\lam_i}{r_i} \cap \mathbb{Z}_- = \emptyset$. \quad \endproof
\begin{remark}
Theorem~\ref{th:VS} enables us to treat all the cases considered in the proof uniformly. 
\end{remark}

If $\lam$ and $r$ are given, then the computation of $f$ requires $\Od(n)$ operations. 
The verification thus require $\Od(n)$ operations. 
\subsection{Computable perturbation bound}\label{sec:EBS}  
As mentioned in Section~\ref{sec:BDS}, the problem of enclosing $\Gamma(A)$ is reduced to that of enclosing $\Gamma(P_1),\dots,\Gamma(P_q)$. 
For $j = 1,\dots,q$, moreover, $P_j$ can be written as $P_j = \widehat{\lam}_jI_{p_j} + Q_j$, where $Q_j \M{p_j}{p_j}$ satisfies $\abs{Q_j} \le R_j$. 
If $p_j = 1$, then we can enclose $\Gamma(P_j)$ by executing the interval variants of Algorithms~\ref{alg:G} and \ref{alg:a} with $\itv{\widehat{\lam}_j}{R_j}$ being the input. 
Otherwise, this approach is not possible. 
In order to enclose $\Gamma(P_j)$ when $p_j \ge 2$, we formulate and prove Theorem~\ref{th:BS}, 
which gives an upper bound for $\normp{\Gamma(\widehat{\lam}_j I_{p_j} + Q_j) - \Gamma(\widehat{\lam}_j I_{p_j})}$, where $p \Zpp \cup \{\infty\}$. 
\begin{theorem}\label{th:BS}  
Let $\omega(\alpha)$ be as in Lemma~\ref{lm:IG}, $p \Zpp \cup \{\infty\}$, $\widehat{\lam}_j \C$, $Q_j \M{p_j}{p_j}$ and $R_j \RMp{p_j}{p_j}$. 
Suppose ${\rm Re}(\widehat{\lam}_j) - \normp{R_j} > 0$ and $\abs{Q_j} \le R_j$, and define 
\begin{eqnarray*}
\delta_p := \normp{R_j}(\Gamma\su{1}({\rm Re}(\widehat{\lam}_j) + \normp{R_j}) + \omega({\rm Re}(\widehat{\lam}_j) + \normp{R_j}) + \omega({\rm Re}(\widehat{\lam}_j) - \normp{R_j})). 
\end{eqnarray*}
Then, $\normp{\Gamma(\widehat{\lam}_j I_{p_j} + Q_j) - \Gamma(\widehat{\lam}_j I_{p_j})} < \delta_p$. 
\end{theorem}
\begin{remark}
We can compute a rigorous upper bound for $\Gamma\su{1}({\rm Re}(\widehat{\lam}_j)+\normp{R_j})$ by slightly modifying Algorithms~\ref{alg:G} and \ref{alg:a}. 
\end{remark}
\proof Let $t \in (0,\infty)$. 
It follows from Lemma~\ref{lm:H} (b) and (c) that 
\begin{eqnarray}
t^{\widehat{\lam}_j I_{p_j} + Q_j - I_{p_j}} - t^{\widehat{\lam}_j I_{p_j} - I_{p_j}} &=& e^{\log(t)((\widehat{\lam}_j - 1)I_{p_j} + Q_j)} - e^{\log(t)(\widehat{\lam}_j - 1)I_{p_j}} \nonumber \\
&=& e^{\log(t)(\widehat{\lam}_j - 1)I_{p_j}}(e^{\log(t)Q_j} - I_{p_j}) = t^{\widehat{\lam}_j - 1}(e^{\log(t)Q_j} - e^0).  \label{eq:t_minus} 
\end{eqnarray}
From $\abs{Q_j} \le R_j$ and Lemma~\ref{lm:H} (d), moreover, we have 
\begin{equation}\label{eq:e_minus}
\normp{e^{\log(t)Q_j} - e^0} \le \abs{\log(t)}\normp{Q_j}e^{\abs{\log(t)}\normp{Q_j}} \le \abs{\log(t)}\normp{R_j}e^{\abs{\log(t)}\normp{R_j}}. 
\end{equation}
The inequality ${\rm Re}(\widehat{\lam}_j) - \normp{R_j} > 0$ gives ${\rm Re}(\widehat{\lam}_j) > 0$, so that $\mu(\widehat{\lam}_j I_{p_j}) \subsetneq \Cpp$. 
The assumption $\abs{Q_j} \le R_j$ yields 
$$
\mu(\widehat{\lam}_j I_{p_j} + Q_j) \subseteq \itv{\widehat{\lam}_j}{\rho(Q_j)} \subseteq \itv{\widehat{\lam}_j}{\rho(\abs{Q_j})} \subseteq \itv{\widehat{\lam}_j}{\rho(R_j)} 
\subseteq \itv{\widehat{\lam}_j}{\normp{R_j}}. 
$$
This and ${\rm Re}(\widehat{\lam}_j) - \normp{R_j} > 0$ give $\mu(\widehat{\lam}_j I_{p_j} + Q_j) \subsetneq \Cpp$. 
The relations $\mu(\widehat{\lam}_j I_{p_j}) \subsetneq \Cpp$, $\mu(\widehat{\lam}_j I_{p_j} + Q_j) \subsetneq \Cpp$, (\ref{eq:Gamma}), (\ref{eq:t_minus}) and (\ref{eq:e_minus}) show 
\begin{eqnarray}
\normp{\Gamma(\widehat{\lam}_j I_{p_j} + Q_j) - \Gamma(\widehat{\lam}_j I_{p_j})} 
&=& \Normp{\int_0^\infty e^{-t}(t^{\widehat{\lam}_j I_{p_j} + Q_j - I_{p_j}} - t^{\widehat{\lam}_j I_{p_j} - I_{p_j}})dt} \nonumber \\
&=& \Normp{\int_0^\infty e^{-t}t^{\widehat{\lam}_j-1}(e^{\log(t)Q_j} - e^0)dt} \nonumber \\
&\le& \int_0^\infty \abs{e^{-t}}\abs{t^{\widehat{\lam}_j-1}}\normp{e^{\log(t)Q_j} - e^0}dt \nonumber \\
&\le& \normp{R_j}\int_0^\infty e^{-t}t^{{\rm Re}(\widehat{\lam}_j)-1}\abs{\log(t)}e^{\abs{\log(t)}\normp{R_j}}dt.  \label{eq:g_minus}
\end{eqnarray}
Let $\mathcal{I}_0 := \int_0^1 e^{-t}t^{{\rm Re}(\widehat{\lam}_j)-1}\abs{\log(t)}e^{\abs{\log(t)}\normp{R_j}}dt$ and  
$\mathcal{I}_\infty := \int_1^\infty e^{-t}t^{{\rm Re}(\widehat{\lam}_j)-1}\abs{\log(t)}$ $e^{\abs{\log(t)}\normp{R_j}}dt$. 
Lemma~\ref{lm:IG} yields
\begin{equation}\label{eq:IG-}
\mathcal{I}_0 = -\gamma\su{1}({\rm Re}(\widehat{\lam}_j) - \normp{R_j}) < \omega({\rm Re}(\widehat{\lam}_j) - \normp{R_j}). 
\end{equation}
From $\Gamma\su{1}(\alpha) = \int_0^\infty e^{-t}t^{\alpha-1}\log(t)dt$ for $\alpha \Rpp$, we moreover have 
\begin{eqnarray}
\mathcal{I}_\infty &=& \Gamma\su{1}({\rm Re}(\widehat{\lam}_j) + \normp{R_j}) - \gamma\su{1}({\rm Re}(\widehat{\lam}_j) + \normp{R_j}) \nonumber \\
&<& \Gamma\su{1}({\rm Re}(\widehat{\lam}_j) + \normp{R_j}) + \omega({\rm Re}(\widehat{\lam}_j) + \normp{R_j}). \label{eq:IG+}
\end{eqnarray}
The relations (\ref{eq:g_minus}) to (\ref{eq:IG+}) prove the inequality. \quad \endproof
\begin{remark}
In \cite{CS}, the estimations $\abs{\log(t)} \le \inv{t}$ and $\abs{\log(t)} \le t$ for $t \in (0,1]$ and $t \in [1,\infty)$, respectively, are used. 
By using the derivatives instead of these estimations, Theorem~\ref{th:BS} gives a smaller bound. 
If we use $\abs{\log(t)} \le \inv{t}$, moreover, then the obtained bound will contain an upper bound for $\gamma({\rm Re}(\widehat{\lam}_j) - \normp{R_j} - 1)$, 
and the condition ${\rm Re}(\widehat{\lam}_j) - \normp{R_j} > 1$ will be required for computing the bound. 
Therefore, the use of the derivatives enables us to weaken the condition. 
On the other hand, using these estimations in \cite{CS} is reasonable. 
This is because the purpose of using these estimations in \cite{CS} is to clarify not quantitative but qualitative properties of $\Gamma(A)$. 
\end{remark}

From Theorem~\ref{th:BS}, we immediately obtain Corollary~\ref{cr:BS}. 
\begin{corollary}\label{cr:BS}
Let $\widehat{\lam}_j$, $Q_j$ and $\delta_p$ be as in Theorem~\ref{th:BS}, and $\delta := \min(\delta_1,\delta_\infty)$.
Under the assumptions in Theorem~\ref{th:BS}, $\Gamma(\widehat{\lam}_j I_{p_j} + Q_j) \in \itv{\Gamma(\widehat{\lam}_j) I_{p_j}}{\delta\om{p_j}}$. 
\end{corollary}
\proof Theorem~\ref{th:BS} and 
$\abs{\Gamma(\widehat{\lam}_j I_{p_j} + Q_j) - \Gamma(\widehat{\lam}_j I_{p_j})} \le \normp{\Gamma(\widehat{\lam}_j I_{p_j} + Q_j) - \Gamma(\widehat{\lam}_j I_{p_j})}\om{p_j}$ for $p = 1,\infty$ 
give $\abs{\Gamma(\widehat{\lam}_j I_{p_j} + Q_j) - \Gamma(\widehat{\lam}_j)I_{p_j}} < \delta\om{p_j}$, proving the result.  \quad \endproof

If ${\rm Re}(\widehat{\lam}_j) - \normp{R_j} > 0$, then the assumption ${\rm Re}(\widehat{\lam}_j) > -5/2$ in Algorithm~\ref{alg:a} is satisfied, 
so that we can enclose $\Gamma(\widehat{\lam}_j)I_{p_j}$ via Algorithms~\ref{alg:G} and \ref{alg:a}. 
%
\subsection{ARs of diagonal blocks}\label{sec:ARS}  
If ${\rm Re}(\widehat{\lam}_j) - \normp{R_j} > 0$ cannot be verified, then Theorem~\ref{th:BS} is not applicable.  
To overcome this issue, we apply the matrix AR based on Lemma~\ref{lm:CS}. 
If the assumption in Theorem~\ref{th:VS} is true, then $\mu(\widehat{\lam}_j I_{p_j} +Q_j) \cap \mathbb{Z}_- = \emptyset$. 
This is because $\mu(\widehat{\lam}_j I_{p_j} +Q_j) = \mu(P_j) \subseteq \mu(A)$ and $\mu(A) \cap \mathbb{Z}_- = \emptyset$, where $P_j$ is as in Section~\ref{sec:BDS}. 
In this case, for $m_j \Zpp$, Lemma~\ref{lm:CS} implies 
\begin{eqnarray}
\Gamma(\widehat{\lam}_j I_{p_j} +Q_j) &=& \inv{\left(\prod_{i=0}^{m_j-1}((\widehat{\lam}_j + i)I_{p_j} + Q_j)\right)}\Gamma((\widehat{\lam}_j + m_j)I_{p_j} + Q_j) \nonumber \\ 
&\in& \inv{\left(\prod_{i=0}^{m_j-1}\itv{(\widehat{\lam}_j + i)I_{p_j}}{R_j}\right)}\Gamma((\widehat{\lam}_j + m_j)I_{p_j} + Q_j),  \label{eq:ARS+}
\end{eqnarray}
provided that any matrix contained in $\prod_{i=0}^{m_j-1}\itv{(\widehat{\lam}_j + i)I_{p_j}}{R_j}$ is nonsingular. 
If we appropriately choose $m_j$, then ${\rm Re}(\widehat{\lam}_j) + m_j - \normp{R_j} > 0$ can be verified, so that Theorem~\ref{th:BS} becomes applicable. 
We can verify nonsingularity of the any matrix, and enclose (\ref{eq:ARS+}) by executing a known verification algorithm. 

If ${\rm Re}(\widehat{\lam}_j) + \normp{R_j} \gg 1$, then the term $\Gamma\su{1}({\rm Re}(\widehat{\lam}_j) + \normp{R_j})$ becomes extremely large. 
In order not to use the large term, we can again execute the AR 
\begin{eqnarray}
\Gamma(\widehat{\lam}_j I_{p_j} +Q_j) &=& \left(\prod_{i=1}^{m_j}((\widehat{\lam}_j-i)I_{p_j} +Q_j)\right)\Gamma((\widehat{\lam}_j - m_j)I_{p_j} +Q_j) \nonumber \\
&\in& \left(\prod_{i=1}^{m_j}\itv{(\widehat{\lam}_j-i)I_{p_j}}{R_j}\right)\Gamma((\widehat{\lam}_j - m_j)I_{p_j} +Q_j). \label{eq:ARS-}
\end{eqnarray}
If ${\rm Re}(\widehat{\lam}_j) - m_j + \normp{R_j} \in [1,2]$, then $\abs{\Gamma\su{1}({\rm Re}(\widehat{\lam}_j) - m_j+\normp{R_j})}$ is not large. 

In (\ref{eq:ARS+}), we need to compute the product $\prod_{i=0}^{m_j-1}\itv{(\widehat{\lam}_j + i)I_{p_j}}{R_j}$. 
If we directly compute this product, then $\Od(m_jp_j^3)$ operations are required, which is prohibitively large when $m_j$ and $p_j$ are large. 
For enclosing this product with only $\Od(m_jp_j^2)$ operations, we present Theorem~\ref{th:ARS}. 
\begin{theorem}\label{th:ARS}  
Let $m_j \Zpp$, $R_j \RMp{p_j}{p_j}$ and $\rr := [\max_i(R_j)_{i1},\dots,\max_i(R_j)_{ip_j}]$. 
Define ${\sf R}_0,\dots,{\sf R}_{m_j-1} \RMp{p_j}{p_j}$ by ${\sf R}_0 := R_j$ and 
$$
{\sf R}_{k+1} := \abs{\widehat{\lam}_j + k + 1}{\sf R}_k + \Abs{\prod_{i=0}^{k}(\widehat{\lam}_j + i)}R_j + \ov{p_j}\rr {\sf R}_k, \quad k = 0,\dots,m_j-2. 
$$
Then, $\prod_{i=0}^{m_j-1}\itv{(\widehat{\lam}_j + i)I_{p_j}}{R_j} \subseteq \itv{(\prod_{i=0}^{m_j-1}(\widehat{\lam}_j+i)) I_{p_j}}{{\sf R}_{m_j-1}}$. 
\end{theorem}
\proof We prove Theorem~\ref{th:ARS} by induction. 
The result is obvious when $m_j = 1$. 
Suppose $\prod_{i=0}^{\ell}\itv{(\widehat{\lam}_j + i)I_{p_j}}{R_j} \subseteq \itv{(\prod_{i=0}^{\ell}(\widehat{\lam}_j + i)) I_{p_j}}{{\sf R}_\ell}$ for $\ell > 1$. 
Then, (\ref{eq:A}) and $\ov{p_j}\rr \ge R_j$ give
\begin{eqnarray*}
\prod_{i=0}^{\ell+1}\itv{(\widehat{\lam}_j + i)I_{p_j}}{R_j} &\subseteq& \itv{(\widehat{\lam}_j + \ell + 1)I_{p_j}}{R_j}\Itv{\left(\prod_{i=0}^{\ell}(\widehat{\lam}_j+i)\right) I_{p_j}}{{\sf R}_\ell} \\
&\subseteq& \Itv{\left(\prod_{i=0}^{\ell+1}(\widehat{\lam}_j + i)\right)I_{p_j}}{\abs{\widehat{\lam}_j + \ell + 1}{\sf R}_\ell + \Abs{\prod_{i=0}^{\ell}(\widehat{\lam}_j + i)}R_j + R_j{\sf R}_\ell} \\
&\subseteq& \Itv{\left(\prod_{i=0}^{\ell+1}(\widehat{\lam}_j + i)\right)I_{p_j}}{{\sf R}_{\ell+1}}.  \quad \endproof
\end{eqnarray*}
The computation of ${\sf R}_{k+1}$ involves $\Od(p_j^2)$ operations for each $k$. 
Therefore, the computation of $\itv{(\prod_{i=0}^{m_j-1}(\widehat{\lam}_j + i)) I_{p_j}}{{\sf R}_{m_j-1}}$ requires only $\Od(m_jp_j^2)$ operations. 

The reduction (\ref{eq:ARS-}) can be accelerated completely analogously. 
\begin{corollary}\label{cr:ARS}  
Let $m_j$, $R_j$ and $\rr$ be as in Theorem~\ref{th:ARS}. 
Define ${\sf S}_1,\dots,{\sf S}_{m_j} \RMp{p_j}{p_j}$ by ${\sf S}_1 := R_j$ and
$$
{\sf S}_{k+1} := \abs{\widehat{\lam}_j - k - 1}{\sf S}_k + \Abs{\prod_{i=1}^{k}(\widehat{\lam}_j - i)}R_j + \ov{p_j}\rr {\sf S}_k, \quad k = 1,\dots,m_j-1. 
$$
Then, $\prod_{i=1}^{m_j}\itv{(\widehat{\lam}_j - i)I_{p_j}}{R_j} \subseteq \itv{(\prod_{i=1}^{m_j}(\widehat{\lam}_j - i)) I_{p_j}}{{\sf S}_{m_j}}$. 
\end{corollary}
%

In practical execution, we need to choose $m_j$. 
We first consider choosing $m_j$ in (\ref{eq:ARS+}). 
As mentioned above, $m_j$ must satisfy ${\rm Re}(\widehat{\lam}_j) + m_j - \normp{R_j} > 0$. 
If $m_j$ is too large, then $\Gamma\su{1}({\rm Re}(\widehat{\lam}_j) + m_j + \normp{R_j}) \gg 1$. 
If ${\rm Re}(\widehat{\lam}_j) + m_j - \normp{R_j}$ is larger than, but close to 0, then $\omega({\rm Re}(\widehat{\lam}_j) + m_j - \normi{R_j}) \gg 1$. 
Based on these observations, we propose determining $m_j = 1 - \lfloor{\rm Re}(\widehat{\lam}_j) - \normi{R_j}\rfloor$, which assures ${\rm Re}(\widehat{\lam}_j) + m_j - \normi{R_j} \in [1,2)$. 
We can analogously choose $m_j$ in (\ref{eq:ARS-}). 
Specifically, we choose $m_j = \lfloor{\rm Re}(\widehat{\lam}_j) - \normi{R_j}\rfloor -1$, which assures ${\rm Re}(\widehat{\lam}_j) - m_j - \normi{R_j} \in [1,2)$. 

There exists the case where the AR is required even when $p_j = 1$. 
To be specific, we can not execute Algorithm~\ref{alg:a} if ${\rm Re}(\widehat{\lam}_j) - R_j > -5/2$ can not be verified. 
In this case, we execute the AR 
\begin{eqnarray*}
\Gamma(\widehat{\lam}_j + Q_j) &=& \inv{\left(\prod_{i=0}^{m_j-1}(\widehat{\lam}_j + i + Q_j)\right)}\Gamma(\widehat{\lam}_j + m_j + Q_j) \\
&\in&  \inv{\left(\prod_{i=0}^{m_j-1}\itv{\widehat{\lam}_j + i}{R_j}\right)}\Gamma(\widehat{\lam}_j + m_j + Q_j) 
\end{eqnarray*}
in order to make ${\rm Re}(\widehat{\lam}_j) - R_j + m_j$ larger than $-5/2$. 
We determine $m_j$ such that $m_j = -2 - \lfloor{\rm Re}(\widehat{\lam}_j) - R_j\rfloor$, which assures ${\rm Re}(\widehat{\lam}_j) - R_j + m_j \in [-2,-1)$. 
%
\subsection{Overall algorithm}\label{sec:OAS}  
Based on Sections~\ref{sec:BDS} to \ref{sec:ARS}, we propose an algorithm for enclosing $\Gamma(A)$. 
\begin{algorithm}\label{alg:S}
Let $P_j$ and $\bm{W}$ be as in Section~\ref{sec:BDS}, and $\bm{P}_j \IM{p_j}{p_j}$ and $\bm{\Gamma}_j \IM{p_j}{p_j}$ contain $P_j$ and $\Gamma(P_j)$, respectively, for $j = 1,\dots,q$. 
This algorithm computes $\bm{\Gamma} \IM{n}{n}$ such that $\bm{\Gamma} \ni \Gamma(A)$. 
If the algorithm successfully terminated, then $\mu(A) \cap \mathbb{Z}_- = \emptyset$ is moreover proved. 
\begin{description}
\item[Step 1.] Compute $\bm{W}$ and $\bm{P}_j$, $j = 1,\dots,q$ by executing \cite[Algorithm~1]{M14}. 
Note that $r$ in Section~\ref{sec:VS} is also obtained in this process. 
\item[Step 2.] Let $f$ be as in Section~\ref{sec:VS}. 
If $\min_i(\abs{f_i} - r_i) > 0$ cannot be verified, terminate with failure. 
Otherwise, $\mu(A) \cap \mathbb{Z}_- = \emptyset$ is proved. 
\item[Step 3.] Compute $\bm{\Gamma}_j \IM{p_j}{p_j}$ for all $j$ by repeatedly executing Algorithm~\ref{alg:EBS}. 
\item[Step 4.] Compute $\bm{\Gamma}$ by $\bm{\Gamma} = \bm{W}{\rm diag}(\bm{\Gamma}_1,\dots,\bm{\Gamma}_q)\inv{\bm{W}}$. 
Terminate. 
\end{description}
\end{algorithm}
\begin{algorithm}\label{alg:EBS}
This algorithms computes $\bm{\Gamma}_j$ in Algorithm~\ref{alg:S}. 
\begin{description}
\item[Step 1.] If $p_j = 1$, then go to Step~2. 
Otherwise, go to Step~3. 
\item[Step 2.] Compute $\bm{\Gamma}_j$ by executing the interval valiants of Algorithms~\ref{alg:G} and \ref{alg:a}, and AR if necessary. 
Terminate. 
\item[Step 3.] If ${\rm Re}(\widehat{\lam}_j) - \normp{R_j} > 0$ for $p = 1,\infty$ cannot be verified, then go to Step~4. 
Otherwise, go to Step~5. 
\item[Step 4.]  Compute $\bm{\Gamma}_j$ with the AR (\ref{eq:ARS+}). 
Terminate. 
\item[Step 5.] Compute $\bm{\Gamma}_j$ with the AR (\ref{eq:ARS-}). 
Terminate. 
\end{description}
\end{algorithm}
Step~1 in Algorithm~\ref{alg:S} involves $\Od(n^3)$ operations (see \cite[Section~3.4]{M14}). 
Step~4 in Algorithm~\ref{alg:S} also involves $\Od(n^3)$ operations. 
The computational cost of Algorithm~\ref{alg:EBS} is $\Od(p_j^3 + m_jp_j^2)$. 
From this and $\sum_{j=1}^qp_j = n$, Step~3 in Algorithm~\ref{alg:S} requires $\Od(n^3 + \sum_{j=1}^qm_jp_j^2)$ operations. 
Costs of other parts in Algorithm~\ref{alg:S} are negligible. 
Algorithm~\ref{alg:S} thus involves only $\Od(n^3)$ operations if $\sum_{j=1}^qm_jp_j^2$ is $\Od(n^3)$. 
\section{Algorithm based on the NJD}\label{sec:J}
Let $N_j$ and $n_j$ be as in Section~\ref{sec:P}. 
When $A$ is defective or close to defective, the matrix $X$ in Section~\ref{sec:BDS} becomes singular or ill-conditioned, which causes failure of \cite[Algorithm~1]{M14}. 
Even in such situations, we can utilize the NJD $AZ \approx ZJ$, where $Z,J \M{n}{n}$, $Z$ is nonsingular, $J = {\rm diag}(J_1,\dots,J_p)$, 
$J_j = \lam_jI_{n_j} + N_j$, $j = 1,\dots,p$, and $\sum_{j=1}^pn_j = n$. 
We proceed similarly to Section~\ref{sec:S}. 
\subsection{The VBD based on the NJD}\label{sec:BDJ}
Let $q$, $i_1\su{1},\dots,i_{p_1}\su{1},\dots,i_1\su{q},\dots,i_{p_q}\su{q}$, $W_j$, $P_j$ and $\bm{W}$ be as in Section~\ref{sec:BDS}. 
We can execute \cite[Algorithm~3]{M14} utilizing the NJD instead of the numerical block diagonalization in \cite[Section~4]{M14}. 
Then, we can obtain $\widehat{\lam}_j \C$, $r \RVp{q}$, $X_j \M{n}{p_j}$ and $\Delta_j \RMp{n}{p_j}$ such that 
$\itv{X_j}{\Delta_j} \ni W_j$, $\itv{\widehat{\lam}_j}{r_j} \supseteq \mu(P_j)$ and $\bigcup_{j=1}^q\itv{\widehat{\lam}_j}{r_j} \supseteq \mu(A)$. 
As byproducts, this algorithm also gives $R_j,M_j \RMp{p_j}{p_j}$ such that $\itv{\widehat{\lam}_jI_{p_j} + M_j}{R_j} \ni P_j$ and $M_j = {\rm diag}(N_{k_1\su{j}},\dots,N_{k_{s_j}\su{j}})$, 
where $\sum_{\ell = 1}^{s_j}n_{k_\ell\su{j}} = p_j$. 
If verification for nonsingularity of any matrix contained in $\bm{W}$ is succeeded, 
then $\bm{W}$ and ${\rm diag}(\itv{\widehat{\lam}_1I_{p_1} + M_1}{R_1},\dots,\itv{\widehat{\lam}_qI_{p_q} + M_q}{R_q})$ can be regarded as the result of the VBD. 
\subsection{Verification of $\bm{\mu(A) \cap \mathbb{Z}_- = \emptyset}$}\label{sec:VJ}  
Similarly to Section~\ref{sec:VS}, we have 
\begin{corollary}\label{cr:VJ}  
Let $\widehat{\lam} \V{q}$ and $r \RVp{q}$ satisfy $\mu(A) \subseteq \bigcup_{i=1}^q\itv{\widehat{\lam}_i}{r_i}$. 
Define $f \V{q}$ by 
$$
f_i := \max({\rm Re}(\widehat{\lam}_i), \lfloor{\rm Re}(\widehat{\lam}_i)\rfloor - {\rm Re}(\widehat{\lam}_i), {\rm Re}(\widehat{\lam}_i) - \lceil{\rm Re}(\widehat{\lam}_i)\rceil) 
+ {\rm Im}(\widehat{\lam}_i){\sf i}, \quad i = 1,\dots,q. 
$$ 
If $\min_i(\abs{f_i} - r_i) > 0$, then $\mu(A) \cap \mathbb{Z}_- = \emptyset$. 
\end{corollary}
%
%
\subsection{Computable perturbation bound}\label{sec:EBJ}  
The diagonal block $P_j$ can be written as $P_j = \widehat{\lam}_j I_{p_j} + M_j + Q_j$, where $Q_j \M{p_j}{p_j}$ satisfies $\abs{Q_j} \le R_j$. 
We can derive an upper bound for $\normp{\Gamma(\widehat{\lam}_j I_{p_j} + M_j + Q_j) - \Gamma(\widehat{\lam}_j I_{p_j} + M_j)}$ analogously to Section~\ref{sec:EBS}. 
\begin{theorem}\label{th:BJ}  
Let $\omega(\alpha)$ be as in Lemma~\ref{lm:IG}, $p \Zpp \cup \{\infty\}$, and $\widehat{\lam}_j$, $Q_j$, $M_j$ and $R_j$ be as above. 
Suppose ${\rm Re}(\widehat{\lam}_j) - \normp{M_j+R_j} > 0$ and $\abs{Q_j} \le R_j$, and define 
\begin{eqnarray*}
\delta_p &:=& \normp{R_j}(\Gamma\su{1}({\rm Re}(\widehat{\lam}_j) + \normp{M_j + R_j}) + \omega({\rm Re}(\widehat{\lam}_j) + \normp{M_j + R_j} ) \\
&&+ \omega({\rm Re}(\widehat{\lam}_j) - \normp{M_j + R_j})).
\end{eqnarray*}
Then, $\normp{\Gamma(\widehat{\lam}_j I_{p_j} + M_j + Q_j) - \Gamma(\widehat{\lam}_j I_{p_j} + M_j)} < \delta_p$. 
\end{theorem}
%

Theorem~\ref{th:BJ} immediately gives Corollary~\ref{cr:BJ}. 
\begin{corollary}\label{cr:BJ}
Let $\widehat{\lam}_j$, $Q_j$, $M_j$ and $\delta_p$ be as in Theorem~\ref{th:BJ}, and $\delta := \min(\delta_1,\delta_\infty)$. 
Under the assumptions in Theorem~\ref{th:BJ}, it holds that $\Gamma(\widehat{\lam}_j I_{p_j} + M_j + Q_j) \in \itv{\Gamma(\widehat{\lam}_j I_{p_j} + M_j)}{\delta\om{p_j}}$.
\end{corollary}
%

We can enclose $\Gamma(\widehat{\lam}_j I_{p_j} + M_j)$ by executing Algorithms~\ref{alg:G} and \ref{alg:a}, because $\Gamma(\widehat{\lam}_j I_{p_j} + M_j)$ is equal to 
$$
{\rm diag}\left(\left[\begin{array}{ccc}
\dsfrac{\Gamma\su{0}(\widehat{\lam}_j)}{0!} & \cdots & \dsfrac{\Gamma\su{n_{k_1\su{j}}-1}(\widehat{\lam}_j)}{(n_{k_1\su{j}}-1)!} \\
 & \ddots & \vdots \\
 & & \dsfrac{\Gamma\su{0}(\widehat{\lam}_j)}{0!} \\
\end{array}\right],\dots,\left[\begin{array}{ccc}
\dsfrac{\Gamma\su{0}(\widehat{\lam}_j)}{0!} & \cdots & \dsfrac{\Gamma\su{n_{k_{s_j}\su{j}}-1}(\widehat{\lam}_j)}{(n_{k_{s_j}\su{j}}-1)!} \\
 & \ddots & \vdots \\
 & & \dsfrac{\Gamma\su{0}(\widehat{\lam}_j)}{0!} \\
\end{array}\right]\right). 
$$
\subsection{ARs of diagonal blocks}\label{sec:ARJ}  
Suppose the assumption in Corollary~\ref{cr:VJ} is true. 
If $M_j = 0$, i.e., $n_{k_1\su{j}} = \cdots = n_{k_{s_j}\su{j}} = 1$, then the ARs in Section~\ref{sec:ARS} are possible. 
Otherwise, we execute the ARs as follows: 
Let $m_j \Zp$. 
If ${\rm Re}(\widehat{\lam}_j) - \normp{M_j+R_j} > 0$ can not be verified, then we execute the AR
\begin{eqnarray}
\Gamma(\widehat{\lam}_j I_{p_j} + M_j + Q_j) &=& \inv{\left(\prod_{i=0}^{m_j-1}((\widehat{\lam}_j + i)I_{p_j} + M_j + Q_j)\right)}\Gamma((\widehat{\lam}_j + m_j)I_{p_j} + M_j + Q_j) \nonumber \\ 
&\in& \inv{\left(\prod_{i=0}^{m_j-1}\itv{(\widehat{\lam}_j + i)I_{p_j} + M_j}{R_j}\right)}\Gamma((\widehat{\lam}_j + m_j)I_{p_j} + M_j + Q_j),  \label{eq:ARJ+}
\end{eqnarray}
provided that any matrix contained in $\prod_{i=0}^{m_j-1}\itv{(\widehat{\lam}_j + i)I_{p_j} + M_j}{R_j}$ is nonsingular. 
If ${\rm Re}(\widehat{\lam}_j)+\normp{M_j+R_j} \gg 1$, alternatively, then we execute 
\begin{eqnarray}
\Gamma(\widehat{\lam}_j I_{p_j} + M_j + Q_j) &=& \left(\prod_{i=1}^{m_j}((\widehat{\lam}_j-i)I_{p_j} + M_j + Q_j)\right)\Gamma((\widehat{\lam}_j - m_j)I_{p_j} + M_j + Q_j) \nonumber \\
&\in& \left(\prod_{i=1}^{m_j}\itv{(\widehat{\lam}_j-i)I_{p_j} + M_j}{R_j}\right)\Gamma((\widehat{\lam}_j - m_j)I_{p_j} + M_j + Q_j).  \label{eq:ARJ-}
\end{eqnarray}

The theories for verifying $\mu(A) \cap \mathbb{Z}_- = \emptyset$ and enclosing $\Gamma(\widehat{\lam}_j I_{p_j} + M_j + Q_j)$ seems to be analogues of those in Section~\ref{sec:S}. 
However, theories for accelerating the ARs are different. 
\begin{theorem}\label{th:ARJ}  
Let $\widehat{\lam}_j$, $R_j$, $M_j$ and $n_{k_\ell\su{j}}$, $\ell = 1,\dots,s_j$ be as above, $\rr$ be as in Theorem~\ref{th:ARS}, $m_j \Zpp$, and $n_{\rm max}\su{j} := \max_\ell n_{k_\ell\su{j}}$. 
For $k = 0,\dots,m_j-2$, define $\beta_0\su{k},\dots,\beta_{n_{\rm max}\su{j}-1}\su{k} \C$ by 
$\beta_0\su{0} := \widehat{\lam}_j$, $\beta_1\su{0} := 1$, $\beta_\ell\su{0} := 0$, $\ell = 2,\dots,n_{\rm max}\su{j}-1$, 
$$
\beta_0^{(k+1)} := (\widehat{\lam}_j+k+1)\beta_0\su{k}, \quad \beta_{\ell}^{(k+1)} := (\widehat{\lam}_j+k+1)\beta_\ell\su{k} + \beta_{\ell-1}\su{k}, \quad \ell = 1,\dots, n_{\rm max}\su{j}-1.
$$
For $\ell = 1,\dots,s_j$ and $k = 0,\dots,m_j-2$, let 
\begin{eqnarray*} 
R_j\su{\ell} &:=& \left[(R_j)_{:\sum_{i=1}^{\ell-1}n_{k_i\su{j}}+1},\dots,(R_j)_{:\sum_{i=1}^{\ell}n_{k_i\su{j}}}\right] \RMp{p_j}{n_{k_\ell\su{j}}},  \\
\rc\su{j,\ell} &:=& \left[\max_{i}(R_j\su{\ell})_{1i},\dots,\max_i(R_j\su{\ell})_{p_ji}\right]^T \RVp{p_j},  \\
w\su{k,\ell} &:=& \left[\abs{\beta_0\su{k}},\abs{\beta_0\su{k}}+\abs{\beta_1\su{k}},\dots,\sum_{i=0}^{n_{k_\ell\su{j}}-1}\abs{\beta_i\su{k}}\right] \RMp{1}{n_{k_\ell\su{j}}}, \\
{\sf Q}_k &:=& \left[\rc\su{j,1}w\su{k,1},\dots,\rc\su{j,s_j}w\su{k,s_j}\right] \RMp{p_j}{p_j}. 
\end{eqnarray*}
Let also $C_k := \sum_{i=0}^{n_{\rm max}\su{j}-1}\beta_i\su{k}M_j^i$ for $k = 0,\dots,m_j-1$. 
Define ${\sf R}_0,\dots,{\sf R}_{m_j-1} \RMp{p_j}{p_j}$ by  ${\sf R}_0 := R_j$ and 
$$
{\sf R}_{k+1} := {\sf Q}_k + \abs{\widehat{\lam}_j+k+1}{\sf R}_k + M_j{\sf R}_k + \ov{p_j}\rr {\sf R}_k, \quad k = 0,\dots,m_j-2. 
$$
Then, $\prod_{i=0}^{m_j-1} \itv{(\widehat{\lam}_j + i)I_{p_j} + M_j}{R_j} \subseteq \itv{C_{m_j-1}}{{\sf R}_{m_j-1}}$. 
\end{theorem}
\proof We prove Theorem~\ref{th:ARJ} by induction. 
The result is obvious when $m_j = 1$. 
Suppose $\prod_{i=0}^\ell\itv{(\widehat{\lam}_j + i)I_{p_j} + M_j}{R_j} \subseteq \itv{C_\ell}{{\sf R}_\ell}$ for $\ell > 1$. 
Then, (\ref{eq:A}) gives
\begin{eqnarray*}
&&\prod_{i=0}^{\ell+1}\itv{(\widehat{\lam}_j + i)I_{p_j} + M_j}{R_j} \\
&&\subseteq \itv{(\widehat{\lam}_j + \ell + 1)I_{p_j} + M_j}{R_j}\itv{C_\ell}{{\sf R}_\ell} \\
&&\subseteq \itv{((\widehat{\lam}_j + \ell + 1)I_{p_j} + M_j)C_\ell}{\abs{(\widehat{\lam}_j + \ell + 1)I_{p_j} + M_j}{\sf R}_\ell+R_j\abs{C_\ell}+R_j{\sf R}_\ell} \\
&&=: \itv{\widehat{C}_{\ell+1}}{\widehat{{\sf R}}_{\ell+1}}. 
\end{eqnarray*}
Since $M_j^{n_{\rm max}\su{j}} = 0$, it follows that 
$$
\widehat{C}_{\ell+1} = ((\widehat{\lam}_j + \ell + 1)I_{p_j} + M_j)\left(\sum_{i=0}^{n_{\rm max}\su{j}-1}\beta_i^{(\ell)}M_j^i\right) = \sum_{i=0}^{n_{\rm max}\su{j}-1}\beta_i^{(\ell+1)}M_j^i  = C_{\ell+1}. 
$$
The term $R_j\abs{C_\ell}$ in $\widehat{{\sf R}}_{\ell+1}$ can be estimated as follows: 
\begin{eqnarray*}
R_j\abs{C_\ell} &=& R_j\sum_{i=0}^{n_{\rm max}\su{j}-1}\abs{\beta_i\su{\ell}}M_j^i  \\
&=& \left[R_j\su{1}\sum_{i=0}^{n_{k_1\su{j}}-1}\abs{\beta_i\su{\ell}}N_{k_1\su{j}}^i,\dots,R_j\su{s_j}\sum_{i=0}^{n_{k_{s_j}\su{j}}-1}\abs{\beta_i\su{\ell}}N_{k_{s_j}\su{j}}^i\right] \\
&\le& \left[\rc\su{j,1}(\ov{n_{k_1\su{j}}})^T\sum_{i=0}^{n_{k_1\su{j}}-1}\abs{\beta_i\su{\ell}}N_{k_1\su{j}}^i,\dots,
\rc\su{j,s_j}(\ov{n_{k_{s_j}\su{j}}})^T\sum_{i=0}^{n_{k_{s_j}\su{j}}-1}\abs{\beta_i\su{\ell}}N_{k_{s_j}\su{j}}^i\right] \\
&=& \left[\rc\su{j,1}w\su{\ell,1},\dots,\rc\su{j,s_j}w\su{\ell,s_j}\right] = {\sf Q}_\ell.  
\end{eqnarray*}
From this and $\ov{p_j}\rr \ge R_j$, we obtain 
$$
\widehat{{\sf R}}_{\ell+1} \le {\sf Q}_\ell + \abs{\widehat{\lam}_j + \ell + 1}{\sf R}_\ell + M_j{\sf R}_\ell + \ov{p_j}\rr{\sf R}_\ell = {\sf R}_{\ell+1}. 
$$
Hence, $\prod_{i=0}^{\ell+1}\itv{(\widehat{\lam}_j + i)I_{p_j} + M_j}{R_j} \subseteq \itv{C_{\ell+1}}{{\sf R}_{\ell+1}}$. \quad \endproof

It is obvious that we do not need to execute the matrix multiplications $M_j^i$ and $M_j{\sf R}_k$ in $C_{m_j-1}$ and ${\sf R}_{k+1}$, respectively, via floating point arithmetic. 
In fact, $C_{m_j-1}$ and $M_j{\sf R}_k$ can be written as follows:
\begin{eqnarray*}
C_{m_j-1} &=& {\rm diag}\left(\left[\begin{array}{ccc}
\beta_0\su{m_j-1} & \cdots & \beta_{n_{k_1\su{j}}-1}\su{m_j-1} \\
 & \ddots & \vdots \\
 & & \beta_0\su{m_j-1}
\end{array}\right], \dots, \left[\begin{array}{ccc}
\beta_0\su{m_j-1} & \cdots & \beta_{n_{k_{s_j}\su{j}}-1}\su{m_j-1} \\
 & \ddots & \vdots \\
 & & \beta_0\su{m_j-1}
\end{array}\right]\right),  \\
M_j{\sf R}_k &=& \left[({\sf R}_k\su{1})^T,\dots,({\sf R}_k\su{s_j})^T\right]^T, \  \mbox{where} \\
{\sf R}_k\su{\ell} &:=& \left[\begin{array}{ccc}
({\sf R}_k)_{{\sum_{i=1}^{\ell-1}n_{k_i\su{j}}+2} \ 1} & \dots & ({\sf R}_k)_{{\sum_{i=1}^{\ell-1}n_{k_i\su{j}}+2} \ p_j} \\
\vdots & \ddots & \vdots \\
({\sf R}_k)_{{\sum_{i=1}^{\ell}n_{k_i\su{j}}} \ 1} & \dots & ({\sf R}_k)_{{\sum_{i=1}^{\ell}n_{k_i\su{j}}} \ p_j} \\
0 & \cdots & 0 \\
\end{array}\right], \quad \ell = 1,\dots,s_j. 
\end{eqnarray*}
Hence, the computations of $C_{m_j-1}$ and ${\sf R}_{m_j-1}$ require only $\Od(m_jp_j^2)$ operations. 

The enclosure of $\prod_{i=1}^{m_j}\itv{(\widehat{\lam}_j-i)I_{p_j} + M_j}{R_j}$ in (\ref{eq:ARJ-}) can also be accelerated. 
\begin{corollary}\label{cr:ARJ}
Let $\widehat{\lam}_j$, $R_j$, $M_j$, $m_j$, $n_{\rm max}\su{j}$, $\rr$ and ${\sf Q}_k$ be as in Theorem~\ref{th:ARJ}. 
For $k = 1,\dots,m_j-1$, define $\beta_0\su{k},\dots,\beta_{n_{\rm max}\su{j}-1}\su{k} \C$ by 
$\beta_0\su{1} := \widehat{\lam}_j - 1$, $\beta_1\su{1} := 1$, $\beta_\ell\su{1} := 0$, $\ell = 2,\dots,n_{\rm max}\su{j}-1$, 
$$
\beta_0^{(k+1)} := (\widehat{\lam}_j-k-1)\beta_0\su{k}, \quad \beta_{\ell}^{(k+1)} := (\widehat{\lam}_j-k-1)\beta_\ell\su{k} + \beta_{\ell-1}\su{k}, \quad \ell = 1,\dots, n_{\rm max}\su{j}-1.
$$
Let $D_k := \sum_{i=0}^{n_{\rm max}\su{j}-1}\beta_i\su{k}M_j^i$ for $k = 1,\dots,m_j$. 
Define ${\sf S}_1,\dots,{\sf S}_{m_j} \RMp{p_j}{p_j}$ by ${\sf S}_1 := R_j$ and 
$$
{\sf S}_{k+1} := {\sf Q}_k + \abs{\widehat{\lam}_j-k-1}{\sf S}_k + M_j{\sf S}_k + \ov{p_j}\rr {\sf S}_k, \quad k = 1,\dots,m_j-1. 
$$
Then, $\prod_{i=1}^{m_j} \itv{(\widehat{\lam}_j - i)I_{p_j} + M_j}{R_j} \subseteq \itv{D_{m_j}}{{\sf S}_{m_j}}$. 
\end{corollary}
%

We can determine $m_j$ in (\ref{eq:ARJ+}) and (\ref{eq:ARJ-}) analogously to Section~\ref{sec:ARS}, where $\normi{R_j}$ is replaced by $\normi{M_j + R_j}$. 
\subsection{Overall algorithm}\label{sec:OAJ}  
Based on Sections~\ref{sec:BDJ} to \ref{sec:ARJ}, we propose Algorithm~\ref{alg:J}. 
\begin{algorithm}\label{alg:J}
Let $P_j$, $\bm{W}$, $\bm{P}_j$, $\bm{\Gamma}_j$ and $\bm{\Gamma}$ be as in Algorithm~\ref{alg:S}. 
This algorithm computes $\bm{\Gamma}$. 
Moreover, $\mu(A) \cap \mathbb{Z}_- = \emptyset$ is proved if successful. 
\begin{description}
\item[Step 1.] Compute $\bm{W}$ and $\bm{P}_j$, $j = 1,\dots,q$ by executing the Jordan valiant of \cite[Algorithm~3]{M14}. 
Note that $r$ in Section~\ref{sec:BDJ} is also obtained. 
\item[Step 2.] Analogous to that in Algorithm~\ref{alg:S}, where $f$ in Section~\ref{sec:VJ} is used instead. 
\item[Step 3.] Compute $\bm{\Gamma}_j$ for all $j$ via Algorithm~\ref{alg:EBJ}. 
\item[Step 4.] Similar to that in Algorithm~\ref{alg:S}. 
\end{description}
\end{algorithm}
\begin{algorithm}\label{alg:EBJ}
This algorithms computes $\bm{\Gamma}_j$ in Algorithm~\ref{alg:J}. 
\begin{description}
\item[Steps 1 and 2.] Similar to those in Algorithm~\ref{alg:EBS}. 
\item[Step 3.] Analogous to that in Algorithm~\ref{alg:EBS}, where $\normp{R_j}$ is replaced by $\normp{M_j + R_j}$. 
\item[Step 4.] Compute $\bm{\Gamma}_j$ with the AR (\ref{eq:ARS+}) if $M_j = 0$, (\ref{eq:ARJ+}) otherwise. 
Terminate. 
\item[Step 5.] Compute $\bm{\Gamma}_j$ with the AR (\ref{eq:ARS-}) if $M_j = 0$, (\ref{eq:ARJ-}) otherwise. 
Terminate. 
\end{description}
\end{algorithm}
The NJD involves $\Od(n^4)$ operations. 
Algorithm~\ref{alg:J} thus involves $\Od(n^4)$ operations if $\sum_{j=1}^qm_jp_j^2$ is $\Od(n^4)$. 
\section{Numerical results}\label{sec:N}
We used a computer with an Intel Core 1.51 GHz CPU, 16.0 GB RAM, and MATLAB R2012a with the Intel Math Kernel Library and IEEE 754 double precision. 
We denote the compared algorithms as follows:
\begin{description}
\item[{\tt  Gs}:] Algorithm~\ref{alg:S}, where $\mu(A) \cap \mathbb{Z}_- = \emptyset$ is verified, 
\item[{\tt  Gj}:] Algorithm~\ref{alg:J}, where $\mu(A) \cap \mathbb{Z}_- = \emptyset$ is verified, and 
\item[{\tt V}:] VERSOFT routine {\tt VERMATFUN}, where $\mu(A) \cap \mathbb{Z}_- = \emptyset$ is not verified.  
\end{description}
In {\tt Gs} and {\tt Gj}, we perform the NSD and NJD by MATLAB and NAClab \cite{Z} routines {\tt eig} and {\tt NumericalJordanForm}, respectively. 
The routine {\tt NumericalJordanForm} generally returns not $Z$ and $J$ but $\widehat{Z}$ and $\widehat{J}$ such that 
$\widehat{A}\widehat{Z} \approx \widehat{Z}\widehat{J}$, $\widehat{J} = {\rm diag}(\widehat{J}_1,\dots,\widehat{J}_p)$, and superdiagonal entries of $\widehat{J}_k$ are not necessarily one. 
However, we can compute $Z$ and $J$ from $\widehat{Z}$ and $\widehat{J}$ (see \cite[Section~5]{M18ex}). 
In {\tt Gs} and {\tt Gj}, we computed products of an interval matrix and an interval matrix containing inverse matrices via {\tt verifylss}.  
In {\tt V}, we called {\tt vermatfun('gamma(z)',A)} when $A$ is Hermitian, invoking the INTLAB routine {\tt gamma}. 
When $A$ is not Hermitian, we called {\tt vermatfun('verGamma(z)',A)}, 
where {\tt verGamma} is a routine which computes an interval containing $\Gamma(z)$ for $z \C$ based on Section~\ref{sec:GD}. 
See \url{http://web.cc.iwate-u.ac.jp/~miyajima/MGF.zip} for details of the implementations, 
where INTLAB codes for {\tt Gs}, {\tt Gj}, {\tt V}, and {\tt verGamma} (denoted by {\tt Gs.m}, {\tt Gj.m}, {\tt V.m}, and {\tt verGamma.m}) are uploaded. 

Let $\itv{\ap{\Gamma}}{R} \ni \Gamma(A)$. 
To assess quality of enclosure, define the relative radius {\rm RR} by ${\rm RR} := \normi{R}/\normi{\ap{\Gamma}}$. 
For some problems, {\tt Gs} or {\tt V} failed. 
The reason for the failure of {\tt Gs} is that \cite[Algorithm~1]{M14} failed because the nonsingularity of $X$ cannot be verified. 
That of {\tt V} is enclosing all the eigenvalues and eigenvectors of $A$ failed. 
\subsection*{Example~1}
We applied the algorithms to four classes of matrices, ``frank'', ``gcdmat'', ``minij'', and ``poisson'', 
available from the MATLAB {\tt gallery} function, and chose matrices of various $n$ for each of the classes. 
For the ``gcdmat'' and ``minij'' matrices, we divided the generated matrix by $n$ in order to avoid overflow. 
Tables \ref{tb:frank} to \ref{tb:poisson} display the RR and CPU times (sec) of the algorithms. 
We see that {\tt Gs} and {\tt Gj} were faster than {\tt V}  in many cases. 
\begin{table}[h]
\begin{center} 
\caption{The RR (left half) and CPU times (sec) (right half) for the ``frank'' matrix.}
\label{tb:frank}
\begin{tabular}{l|lll|lll} 
\hline
$n$ & {\tt Gs} & {\tt Gj} & {\tt V} & {\tt Gs} & {\tt Gj} & {\tt V} \\
\hline
5 & 5.9e--12 & 3.8e--12 & 3.6e--12 & 6.8e--1 & 9.0e--1 & 7.6e--1 \\ 
7 & 2.4e--9 & 4.4e--11 & 2.4e--11 & 9.5e--1 & 1.0e+0 & 1.0e+0 \\ 
9 & 2.7e--6 & 2.4e--6 & 7.3e--8 & 9.1e--1 & 1.1e+0 & 1.3e+0 \\ 
11 & 1.4e+0 & 1.8e--2 & 6.0e--4 & 9.5e--1 & 1.2e+0 & 1.5e+0 \\ 
\hline
\end{tabular}
\end{center}
\end{table}
\begin{table}[h]
\begin{center} 
\caption{The RR (left half) and CPU times (sec) (right half) for the ``gcdmat'' matrix divided by $n$.}
\label{tb:gcdmat}
\begin{tabular}{l|lll|lll} 
\hline
$n$ & {\tt Gs} & {\tt Gj} & {\tt V} & {\tt Gs} & {\tt Gj} & {\tt V} \\
\hline
100 & 7.5e--12 & 3.6e--11 & 5.5e--11 & 6.4e--1 & 1.7e+0 & 8.7e+0 \\ 
200 & 3.2e--11 & 2.6e--11 & 1.1e--10 & 2.3e+0 & 6.5e+0 & 3.3e+1 \\ 
300 & 1.1e--10 & 3.3e--11 & 3.2e--10 & 8.0e+0 & 2.5e+1 & 7.4e+1 \\ 
400 & 2.2e--10 & 1.3e--10 & 7.4e--10 & 3.8e+1 & 1.1e+2 & 1.6e+2 \\ 
\hline
\end{tabular}
\end{center}
\end{table}
\begin{table}[h]
\begin{center} 
\caption{The RR (left half) and CPU times (sec) (right half) for the ``minij'' matrix divided by $n$.}
\label{tb:minij}
\begin{tabular}{l|lll|lll} 
\hline
$n$ & {\tt Gs} & {\tt Gj} & {\tt V} & {\tt Gs} & {\tt Gj} & {\tt V} \\
\hline
100 & 1.0e--9 & 9.5e--10 & 1.2e--8 & 7.0e--1 & 1.7e+0 & 1.0e+1 \\ 
200 & 8.1e--9 & 1.9e--7 & 2.2e--7 & 2.6e+0 & 7.8e+0 & 3.8e+1 \\ 
300 & 3.7e--8 & 5.0e--8 & 1.2e--6 & 9.5e+0 & 3.0e+1 & 8.8e+1 \\ 
400 & 8.6e--8 & 1.4e--7 & 4.2e--6 & 3.8e+1 & 1.1e+2 & 1.5e+2 \\ 
\hline
\end{tabular}
\end{center}
\end{table}
\begin{table}[h]
\begin{center} 
\caption{The RR (left half) and CPU times (sec) (right half) for the ``poisson'' matrix.}
\label{tb:poisson}
\begin{tabular}{l|lll|lll} 
\hline
$n$ & {\tt Gs} & {\tt Gj} & {\tt V} & {\tt Gs} & {\tt Gj} & {\tt V} \\
\hline
9 & 2.5e--14 & 2.8e--14 & 1.0e--1 & 1.7e+0 & 2.5e+0 & 1.1e--1 \\ 
36 & 3.7e--13 & 5.0e--13 & failed & 8.3e+0 & 1.2e+1 & failed \\ 
81 & 1.9e--12 & 2.0e--12 & failed & 2.3e+1 & 2.9e+1 & failed \\ 
144 & 8.6e--12 & 1.4e--11 & failed & 4.3e+1 & 5.3e+1 & failed \\ 
\hline
\end{tabular}
\end{center}
\end{table}
\subsection*{Example~2}
We consider the case where $A$ comes close to being defective. 
We applied the algorithms to the problem in \cite[Experiment~1]{FHI}, in which 
$$
A = \left[\begin{array}{cc}
1 & 1 \\
0 & 1 + \varepsilon \\
\end{array}\right], \quad \mbox{where} \quad \varepsilon \ge 0, 
$$
whose eigenvector matrix becomes increasingly ill-conditioned as $\varepsilon \to 0$. 
Table~\ref{tb:ep} reports quantities similar to those in Tables~\ref{tb:frank} to \ref{tb:poisson} with $\varepsilon$ varying from $2^0$ to $2^{-52}$. 
This table shows that the RR by {\tt Gj} stayed about the same, whereas those by {\tt Gs} and {\tt V} increased as $\varepsilon$ decreased. 
\begin{table}[h]
\begin{center} 
\caption{The RR (left half) and CPU times (sec) (right half) in Example~2.}
\label{tb:ep}
\begin{tabular}{l|lll|lll} 
\hline
$\varepsilon$ & {\tt Gs} & {\tt Gj} & {\tt V} & {\tt Gs} & {\tt Gj} & {\tt V} \\
\hline
$2^{0}$ & 2.6e--13 & 2.6e--13 & 2.6e--13 & 3.4e--1 & 3.2e--1 & 3.3e--1 \\ 
$2^{-26}$ & 4.2e--6 & 3.2e--13 & 4.2e--6 & 3.4e--1 & 6.1e--1 & 3.8e--1 \\ 
$2^{-39}$ & 3.4e--2 & 2.9e--13 & 3.4e--2 & 3.4e--1 & 6.0e--1 & 3.8e--1 \\ 
$2^{-48}$ & 1.5e+1 & 3.0e--13 & 1.6e+1 & 3.8e--1 & 6.1e--1 & 3.4e--1 \\ 
$2^{-52}$ & 1.3e+3 & 9.7e--13 & failed & 6.5e--1 & 6.7e--1 & failed \\ 
\hline
\end{tabular}
\end{center}
\end{table}
\subsection*{Example~3}
We consider the case where $A$ is defective. 
Let 
$$A_0 := \left[\begin{array}{rrrr}
 2& 2&1&0\\
 0& 1&1&1\\
-1&-1&0&0\\
 1& 1&1&1
\end{array}\right], ~{\rm whose~Jordan~canonical~form~is}~
\left[\begin{array}{rrrr}
1&1&0&0\\
0&1&1&0\\
0&0&1&1\\
0&0&0&1
\end{array}\right].
$$
We set $A = \sigma A_0$ for a parameter $\sigma \R$. 
Table~\ref{tb:defective} reports quantities similar to those in Table~\ref{tb:ep} for various $\sigma$, showing that {\tt Gj} succeeded for all the problems. 
\begin{table}[h]
\begin{center} 
\caption{The RR (left half) and CPU times (sec) (right half) in Example~3.}
\label{tb:defective}
\begin{tabular}{l|lll|lll} 
\hline
$\sigma$ & {\tt Gs} & {\tt Gj} & {\tt V} & {\tt Gs} & {\tt Gj} & {\tt V} \\
\hline
$2^{-1}$ & failed & 4.0e--12 & failed & failed & 6.4e--1 & failed \\ 
$2^{0}$ & failed & 1.0e--11 & failed & failed & 6.2e--1 & failed \\ 
$2^{1}$ & failed & 1.3e--12 & failed & failed & 6.3e--1 & failed \\ 
$2^{2}$ & failed & 1.2e--12 & failed & failed & 6.4e--1 & failed \\ 
$2^{3}$ & failed & 1.1e--12 & failed & failed & 6.5e--1 & failed \\ 
\hline
\end{tabular}
\end{center}
\end{table}
\subsection*{Example~4}
Consider the case where $A$ is derogatory.  
Let $v\su{j} := (I_8)_{:j}$, $j = 1,\dots,8$, and $P := [v^{(7)},v^{(5)},v^{(3)},v^{(1)},v^{(8)},v^{(6)},v^{(4)},v^{(2)}] \R^{8 \times 8}$. 
Then, $P$ is orthogonal. 
Using $A_0$ and $\sigma$ in Example~3, we set $A = \sigma P{\rm diag}(A_0,A_0)P^T$. 
Table~\ref{tb:derogatory} displays quantities similar to those in Table~\ref{tb:defective}, which also shows the robustness of {\tt Gj}. 
\begin{table}[h]
\begin{center} 
\caption{The RR (left half) and CPU times (sec) (right half) in Example~4.}
\label{tb:derogatory}
\begin{tabular}{l|lll|lll} 
\hline
$\sigma$ & {\tt Gs} & {\tt Gj} & {\tt V} & {\tt Gs} & {\tt Gj} & {\tt V} \\
\hline
$2^{-1}$ & failed & 1.7e--11 & failed & failed & 1.1e+0 & failed \\ 
$2^{0}$ & failed & 2.1e--11 & failed & failed & 1.0e+0 & failed \\ 
$2^{1}$ & failed & 2.6e--12 & failed & failed & 8.9e--1 & failed \\ 
$2^{2}$ & failed & 7.0e--12 & failed & failed & 9.0e--1 & failed \\ 
$2^{3}$ & failed & 9.5e--12 & failed & failed & 8.9e--1 & failed \\ 
\hline
\end{tabular}
\end{center}
\end{table}
\section{Concluding remarks}\label{sec:C}
We have established the new framework for enclosing matrix functions based on the VBD, proposed Algorithms~\ref{alg:S} and \ref{alg:J}, and reported the numerical results. 
As mentioned in Section~\ref{sec:I}, these algorithms are first ones which encloses a matrix function based on this framework. 
Let $\varphi: \mathbb{C}\to\mathbb{C}$ be defined on $\mu(A)$. 
Essentially, we can enclose $\varphi(A)$ based on the VBD framework if the followings are possible: 
\begin{itemize}
\item enclosing $\varphi\su{0}(z)/0!,\dots,\varphi\su{\ell}(z)/\ell!$ for $z \C$ and $\ell \Zp$, and 
\item computing rigorous upper bounds for $\normp{\varphi(\widehat{\lam}_j I_{p_j} + Q_j) - \varphi(\widehat{\lam}_j I_{p_j})}$ and/or 
$\normp{\varphi(\widehat{\lam}_j I_{p_j} + M_j + Q_j) - \varphi(\widehat{\lam}_j I_{p_j} + M_j)}$ for $p = 1,\infty$.   
\end{itemize}
Since $\widehat{\lam}_j I_{p_j}$ and $\widehat{\lam}_j I_{p_j} + M_j$ have simple structures, the derivations of the bounds are easier than those for general matrices. 
For example, enclosing $e^A$, $\sin A$ and $\cos A$ will be possible based on this framework. 
Our future work will be to develop algorithms for enclosing the matrix beta and Bessel functions. 
%
%


\begin{thebibliography}{10}
\bibitem{A}
{\sc H. Arndt}, 
{\em On the interval systems $[x] = [A][x] + [b]$ and the powers of interval matrices in complex interval arithmetics}, 
Reliab. Comput., 13 (2007), pp. 245--259. 
%
\bibitem{CS}
{\sc J.R. Cardoso and A. Sadeghi}, 
{\em Computation of matrix gamma function}, 
BIT, 59 (2019), pp. 343--370. 
%
\bibitem{FHI}
{\sc M. Fasi, N.J. Higham, and B. Iannazzo}, 
{\em An algorithm for the matrix Lambert $W$ function}, 
SIAM J. Matrix Anal. Appl., 36 (2015), pp. 669--685. 
%
\bibitem{H08}
{\sc N.J. Higham}, 
{\em Functions of Matrices: Theory and Computation}, 
SIAM Publications, Philadelphia, 2008. 
%
\bibitem{K}
{\sc W. Kr\"amer}, 
{\em Computation of the gamma function $\Gamma(x)$ for real point and interval arguments}, 
Z. Angew. Math. Mech., 70(6) (1990), pp. 581--584. 
%
\bibitem{M14}
{\sc S. Miyajima}, 
{\em Fast enclosure for all eigenvalues and invariant subspaces in generalized eigenvalue problems}, 
SIAM J. Matrix Anal. Appl., 35 (2014), pp. 1205--1225. 
%
\bibitem{M18rt}
\sameauthor, 
{\em Fast verified computation for the matrix principal $p$th root}, 
J. Comput. Appl. Math., 330 (2018), pp. 276--288. 
%
\bibitem{M18ex}
\sameauthor, 
{\em Verified computation of the matrix exponential}, 
Adv. Comput. Math., 45 (2019), pp. 137--152.
%
\bibitem{M19ln}
\sameauthor, 
{\em Verified computation for the matrix principal logarithm}, 
Linear Algebra Appl., 569 (2019), pp. 38--61. 
%
\bibitem{M19w}
\sameauthor, 
{\em  Verified computation for the matrix Lambert $W$ function}, 
Appl. Math. Comput., 362 (2019), 124555.
%
\bibitem{R}
{\sc J. Rohn}, 
{\em VERSOFT: Verification Software in MATLAB/INTLAB}, 
\url{http://uivtx.cs.cas.cz/~rohn/matlab}
%
\bibitem{R99}
{\sc S.M. Rump}, 
{\em INTLAB - INTerval LABoratory}, 
in Developments in Reliable Computing, T. Csendes, ed., Kluwer Academic Publishers, Dordrecht, 1999, pp. 77--107. 
%
\bibitem{R14}
\sameauthor, 
{\em Verified sharp bounds for the real gamma function over the entire floating-point range}, 
NOLTA, IEICE, 5(3) (2014), pp. 339--348.
%
\bibitem{S}
{\sc J. Spouge}, 
{\em Computation of the gamma, digamma, and trigamma functions}, 
SIAM J. Numer. Anal., 31(3) (1994), pp. 931--944. 
%
\bibitem{YOO}
{\sc N. Yamanaka, T. Okayama, and S. Oishi}, 
{\em Verified error bounds for the real gamma function using double exponential formula over semi-infinite interval}, 
Lect. Notes Comput. Sci., 9582 (2016), pp. 224--230.  
%
\bibitem{Z}
{\sc Z. Zeng and T.-Y. Li}, 
{\em NAClab: A Matlab toolbox for numerical algebraic computation}, 
ACM Commun. Comput. Algebra, 47 (2013), pp. 170--173. 
\end{thebibliography}
\end{document}